\documentclass[12pt]{amsart}
\usepackage{amssymb}

\newcommand{\Bgp}{{\Z^\N}}

\long\def\forget#1\forgotten{}
\newcommand{\issuenumber}{23}
\newcommand{\issuemonth}{December}
\newcommand{\issueyear}{2007}

\newtheorem{thm}{Theorem}[section]
\newtheorem{prob}[thm]{Problem}

\newtheorem{issue}{Issue}

\theoremstyle{definition}

\theoremstyle{remark}

\newcommand{\ed}{
\newpage

\section{Unsolved problems from earlier issues}

{\scriptsize

\begin{issue}
Is $\binom{\Omega}{\Gamma}=\binom{\Omega}{\Tau}$?
\end{issue}

\begin{issue}
Is $\ufin(\cO,\Omega)=\sfin(\Gamma,\Omega)$?
And if not, does $\ufin(\cO,\Gamma)$ imply
$\sfin(\Gamma,\Omega)$?
\end{issue}

\stepcounter{issue}

\begin{issue}
Does $\sone(\Omega,\Tau)$ imply $\ufin(\Gamma,\Gamma)$?
\end{issue}

\begin{issue}
Is $\fp=\fp^*$? (See the definition of $\fp^*$ in that issue.)
\end{issue}

\begin{issue}
Does there exist (in ZFC) an uncountable set satisfying $\sfin(\B,\B)$?
\end{issue}

\stepcounter{issue}

\begin{issue}
Does $X \nin \NON(\M)$ and $Y\nin\mathsf{D}$ imply that
$X\cup Y\nin \COF(\M)$?
\end{issue}

\begin{issue}[CH]
Is $\split(\Lambda,\Lambda)$ preserved under finite unions?
\end{issue}

\begin{issue}
Is $\cov(\M)=\fo$? (See the definition of $\fo$ in that issue.)
\end{issue}

\begin{issue}
Does $\sone(\Gamma,\Gamma)$ always contain an element of cardinality $\fb$?
\end{issue}

\begin{issue}
Could there be a Baire metric space $M$ of weight $\aleph_1$ and a partition
$\mathcal{U}$ of $M$ into $\aleph_1$ meager sets where for each ${\mathcal U}'\subset\mathcal U$,
$\bigcup {\mathcal U}'$ has the Baire property in $M$?
\end{issue}

\stepcounter{issue} 

\begin{issue}
Does there exist (in ZFC) a set of reals $X$ of cardinality $\fd$ such that all
finite powers of $X$ have Menger's property $\sfin(\cO,\cO)$?
\end{issue}

\begin{issue}
Can a Borel non-$\sigma$-compact group be generated by a Hurewicz subspace?
\end{issue}

\begin{issue}[MA]
Is there an uncountable $X\sbst\R$ satisfying $\sone(\BO,\BG)$?
\end{issue}

\begin{issue}[CH]
Is there a totally imperfect $X$ satisfying $\ufin(\cO,\Gamma)$
that can be mapped continuously onto $\Cantor$?
\end{issue}

\begin{issue}[CH]
Is there a Hurewicz $X$ such that $X^2$ is Menger but not Hurewicz?
\end{issue}

\begin{issue}
Does the Pytkeev property of $C_p(X)$ imply that $X$ has Menger's property?
\end{issue}

\begin{issue}
Does every hereditarily Hurewicz space satisfy $\sone(\BG,\BG)$?
\end{issue}

\begin{issue}[CH]
Is there a Rothberger-bounded $G\le\Bgp$ such that $G^2$ is not Menger-bounded?
\end{issue}

\begin{issue}
Let $\cW$ be the van der Waerden ideal.
Are $\cW$-ultrafilters closed under products?
\end{issue}

\begin{issue}
Is the $\delta$-property equivalent to the $\gamma$-property $\binom{\Omega}{\Gamma}$?
\end{issue}

} 

\general\end{document}}

\renewcommand{\>}{\right >}
\newcommand{\Cantor}{{\{0,1\}^\N}}

\newcommand{\fb}{\mathfrak{b}}

\newcommand{\fd}{\mathfrak{d}}
\newcommand{\fp}{\mathfrak{p}}

\newcommand{\NON}{{\mathsf   {NON}}}
\newcommand{\COF}{{\mathsf   {COF}}}

\newcommand{\cN}{\mathcal{N}}

\newcommand{\M}{\mathcal{M}}

\newcommand{\cov}{\mathsf{cov}}

\newcommand{\R}{\mathbb{R}}

\newcommand{\fo}{\mathfrak{od}}

\newcommand{\w}{\omega}

\renewcommand{\split}{\mathsf{Split}}
\newcommand{\bq}{\begin{quote}}
\newcommand{\eq}{\end{quote}}
\newcommand{\cO}{\mathcal{O}}
\newcommand{\B}{\mathcal{B}}
\newcommand{\BG}{\B_\Gamma}

\newcommand{\BO}{\B_\Omega}

\newcommand{\sone}{\mathsf{S}_1}    \newcommand{\sfin}{\mathsf{S}_{fin}}

\newcommand{\ufin}{\mathsf{U}_{fin}}
    
\newcommand{\seq}[1]{\{#1\}_{n\in\N}}
\newcommand{\Union}{\bigcup}
\newcommand{\nin}{\not\in}
\newcommand{\cF}{\mathcal{F}}

\newcommand{\cU}{\mathcal{U}}

\newcommand{\cW}{\mathcal{W}}

\newcommand{\NN}{{\N^\N}}
\newcommand{\N}{\mathbb{N}}
\newcommand{\Z}{\mathbb{Z}}

\newcommand{\sbst}{\subseteq}
\newcommand{\by}[2]{\par\hfill\emph{#1}, #2}
\newcommand{\nby}[1]{\par\hfill\emph{#1}}
\newcommand{\Tau}{\mathrm{T}}
\newcommand{\CE}{\textsc{CE}}

\newcommand{\be}{\begin{enumerate}}
\newcommand{\ee}{\end{enumerate}}
\newcommand{\bi}{\begin{itemize}}
\newcommand{\ei}{\end{itemize}}

\newcommand{\general}{\small\vfill\par\noindent\hrulefill\par
\noindent\textbf{Previous issues.} The previous issues of this
bulletin are available online at\\
\texttt{http://front.math.ucdavis.edu/search?\&t=\%22SPM+Bulletin\%22}
\\[0.1cm]
\textbf{Contributions.} Announcements, discussions, and open problems should be emailed
to \texttt{tsaban@math.biu.ac.il}\\[0.1cm]
\textbf{Subscription.}
To receive this bulletin (free) to your e-mailbox, e-mail us.
}

\newcommand{\link}[1]{\par\hfill{\texttt{#1}}}

\newcommand{\arXiv}[5]{\subsection{#2}{#4}\par\hfill{\arx{#1}}\par\hfill\emph{#3}}

\newcommand{\naAMSPaper}[3]{\subsection{#1}~\par\hfill{\texttt{#3}}\par\hfill\emph{#2}}

\newcommand{\arx}[1]{\texttt{http://arxiv.org/abs/#1}}
\newcommand{\url}[1]{\bq\texttt{#1}\eq}
\newcommand{\online}[1]{The paper is available online at \url{#1}}

\title[$\mathcal{SPM}$ Bulletin \textbf{\issuenumber} (\issuemonth{} \issueyear)]{%
$\mathcal{SPM}$ Bulletin\\[0.5cm]
Issue number \issuenumber: \issuemonth{} \issueyear{} \CE{}}

\begin{document}
\maketitle

\tableofcontents

\section{Editor's note}

A surprising number of new results and directions in ``core'' SPM by
Babinkostova and Scheepers in the last quarter of the year!
See \S\S \ref{bs1}--\ref{bs2} and \ref{BM}--\ref{ssc}.

\medskip

Have a good 2008,

\by{Boaz Tsaban}{tsaban@math.biu.ac.il}

\hfill \texttt{http://www.cs.biu.ac.il/\~{}tsaban}

\section{Research announcements}

\subsection{Cardinal invariants of the continuum and combinatorics on
uncountable cardinals}\label{joerg}

We explore the connection between combinatorial principles on
uncountable cardinals, like stick and club, on the one hand, and
the combinatorics of sets of reals and, in particular, cardinal
invariants of the continuum, on the other hand. For example, we
prove that additivity of measure implies that Martin's axiom holds
for any Cohen algebra. We construct a model in which club holds,
yet the covering number of the null ideal $\cov(\cN)$ is large. We show
that for uncountable cardinals $\kappa\le\lambda$ and $\cF \sbst [\lambda]^\kappa$, if all subsets
of $\lambda$ either contain, or are disjoint from, a member of $\cF$, then $\cF$
has size at least $\cov(\cN)$ etc. As an application, we solve the
Gross space problem under $c = \aleph_2$ by showing that there is such a
space over any countable field. In two appendices, we solve
problems of Fuchino, Shelah and Soukup, and of Kraszewski,
respectively.

\nby{J\"org Brendle}

\arXiv{0709.2893}
{Selection principles and countable dimension\label{bs1}}
{Liljana Babinkostova and Marion Scheepers}
{We characterize countable dimensionality and strong countable dimensionality
by means of an infinite game.}

\arXiv{0709.2895}
{Products and selection principles\label{bs2}}
{Liljana Babinkostova and Marion Scheepers}
{The product of a Sierpinski set and a Lusin set has Menger's property. The
product of a gamma set and a Lusin set has Rothberger's property.}

\arXiv{0709.3016}
{On completely donut (doughnut) sets}
{Piotr Kalemba, Szymon Plewik and Anna Wojciechowska}
{A set $\<A,B\> = \{X \in [\omega]^\omega : A \subseteq X \subseteq B \}$ is a
donut, whenever $A \subseteq B\subseteq \omega$ and $B\setminus A $ is
infinite. A subset $S \subseteq [\omega]^\omega$ is completely donut, whenever
for each donut $\<A,B\>$ there exists a donut $\<C,D\>\subseteq \<A,B\>$ such that $
\<C,D\>\subseteq S$ or $\<C,D\> \cap S = \emptyset.$ If always holds $ \<C,D\> \cap S
= \emptyset$, then $S$ is nowhere donut. We examine families of completely
donut and nowhere donut sets. The results correspond to completely Ramsey and
nowhere Ramsey sets.}

\arXiv{0710.0152}
{On minimal non-potentially closed subsets of the plane}
{Dominique Lecomte\footnote{Lecomte has
recently uploaded quite a few interesting papers to the ArXiv.
Check there.}}
{We study the Borel subsets of the plane that can be made closed by refining
the Polish topology on the real line. These sets are called potentially closed.
We first compare Borel subsets of the plane using products of continuous
functions. We show the existence of a perfect antichain made of minimal sets
among non-potentially closed sets. We apply this result to graphs, quasi-orders
and partial orders. We also give a non-potentially closed set minimum for
another notion of comparison. Finally, we show that we cannot have injectivity
in the Kechris-Solecki-Todorcevic dichotomy about analytic graphs.

Topology and its Applications 154 (2007), 241--262.
}

\arXiv{0710.1085}
{One Dimensional Locally Connected $S$-spaces}
{Joan E. Hart Kenneth Kunen}
{We construct, assuming Jensen's principle $\diamondsuit$, a one-dimensional locally
connected hereditarily separable continuum without convergent sequences. The
construction is an inverse limit in $\omega_1$ steps, and is patterned after the
original Fedorchuk construction of a compact $S$-space. To make it
one-dimensional, each space in the inverse limit is a copy of the Menger
sponge.
}

\arXiv{0710.1402}
{Covering an uncountable square by countably many continuous functions}
{Wieslaw Kubis}
{ We prove that there exists a countable family of continuous real functions
whose graphs together with their inverses cover an uncountable square, i.e. a
set of the form $X\times X$, where $X$ is uncountable. This is motivated by an
old result of Sierpi\'nski, saying that $\aleph_1\times\aleph_1$ is covered by
countably many graphs of functions and inverses of functions. Another
motivation comes from Shelah's study of planar Borel sets without perfect
rectangles.
}

\arXiv{0710.2347}
{Ramsey degrees of finite ultrametric spaces, ultrametric Urysohn spa\-ces and dynamics of their isometry groups}
{L. Nguyen Van Th\'e}
{We study Ramsey-theoretic properties of several natural classes of finite
ultrametric spaces, describe the corresponding Urysohn spaces and compute a
dynamical invariant attached to their isometry groups.}

\arXiv{0710.2352}
{Big Ramsey degrees and divisibility in classes of ultrametric spaces}
{L. Nguyen Van Th\'e}
{Given a countable set $S$ of positive reals, we study finite-dimensional
Ramsey-theoretic properties of the countable ultrametric Urysohn space with
distances in $S$.}

\arXiv{0711.0162}
{Definable Davies' Theorem}
{Asger Tornquist and William Weiss}
{We prove the following analogue of a Theorem of R.O. Davies: Every
$\Sigma^1_2$ function $f:\R\times\R\to\R$ can be represented as a sum of
rectangular $\Sigma^1_2$ functions if and only if all reals are constructible.}

\arXiv{0711.1104\label{BM}}
{Selection Principles and Baire spaces}
{Marion Scheepers}
{We prove that if $X$ is a separable metric space with the Hurewicz covering
property, then the Banach-Mazur game played on $X$ is determined.
The implication is not true when ``Hurewicz covering property'' is
replaced with ``Menger covering property''.}

\arXiv{0711.1322}
{Selective screenability in topological groups\label{scgp}}
{Liljana Babinkostova}
{We examine the selective screenability property in topological groups. In the
metrizable case we also give characterizations in terms of the Haver property
and finitary Haver property respectively relative to left-invariant metrics. We
prove theorems stating conditions under which the properties are preserved by
products. Among metrizable groups we characterize the ones of countable
covering dimension by a natural game.}

\arXiv{0711.1516}
{Selective screenability and the Hurewicz property\label{ssc}}
{Liljana Babinkostova}
{We characterize the Hurewicz covering property in metrizable spaces in terms
of properties of the metrics of the space. Then we show that a weak version of
selective screenability, when combined with the Hurewicz property, implies
selective screenability.}

\naAMSPaper{Partitioning triples and partially ordered sets}
{Albin Jones}
{http://www.ams.org/journal-getitem?pii=S0002-9939-07-09170-8}

\naAMSPaper{A polarized partition relation for cardinals of countable cofinality}
{Albin Jones}
{http://www.ams.org/journal-getitem?pii=S0002-9939-07-09143-5}

\subsection{On topological spaces of singular density and minimal weight}
In a recent paper, Juh\'asz and Shelah establish the consistency
of a regular hereditarily Lindel\"of space of density $\aleph_{\omega_1}$.

It is natural to ask what is the minimal possible value for the
weight of such space; more specifically,
can the weight be $\aleph_{\omega_1}$?

In this paper, we isolate a certain consequence of the
Generalized Continuum Hypothesis,
which we will refer to as the
\emph{Prevalent Singular Cardinals Hypothesis},
and show it implies that every topological space
of density and weight
$\aleph_{\omega_1}$ is not hereditarily Lindel\"of.

The assertion \textsf{PSH} is very weak,
and in fact holds in all currently known models of ZFC.
\link{dx.doi.org/10.1016/j.topol.2007.09.013}
\nby{Assaf Rinot}

\arXiv{0711.4400}
{There is a van Douwen MAD family}
{Dilip Raghavan}
{We prove in ZFC that there is a MAD family of functions in $\NN$ which
is also maximal with respect to infinite partial functions. This solves a long
standing question of van Douwen. We also prove that such families cannot be
analytic. This strengthens Steprans' result that strongly MAD families cannot
be analytic.
}

\arXiv{0712.0584}
{Cardinal sequences of LCS spaces under GCH}
{Juan Carlos Martinez, Lajos Soukup}
{We give full characterization of the sequences of regular cardinals that may
arise as cardinal sequences of locally compact scattered spaces under GCH.
The proofs are based on constructions of universal locally compact scattered
spaces.}

\arXiv{0712.2112}
{Dirichlet sets and Erdos-Kunen-Mauldin theorem}
{Peter Elias}
{By a theorem proved by Erdos, Kunen and Mauldin, for any nonempty perfect set
$P$ on the real line there exists a perfect set $M$ of Lebesgue measure zero
such that $P+M=\mathbb{R}$. We prove a stronger version of this theorem in
which the obtained perfect set $M$ is a Dirichlet set. Using this result we
show that for a wide range of families of subsets of the reals, all additive
sets are perfectly meager in transitive sense. We also prove that every proper
analytic subgroup $G$ of the reals is contained in an $F_\sigma$ set $F$ such that
$F+G$ is a meager null set.}

\arXiv{0712.2393}
{Local Ramsey theory: An abstract approach}
{Jos\'e Mijares and Jes\'us Nieto}
{It is shown that the known notion of selective coideal can be extended to a
family $\mathcal{H}$ of subsets of $\mathcal{R}$, where $(\mathcal{R},\leq,r)$
is a topological Ramsey space in the sense of Todorcevic.
Then it is proven that, if $\mathcal{H}$ selective, the $\mathcal{H}$-Ramsey
and $\mathcal{H}$-Baire subsets of $\mathcal{R}$ are equivalent. This extends
results of Farah for semiselective coideals of
$\mathbb{N}$. Also, it is proven that the family of $\mathcal{H}$-Ramsey subsets
of $\mathcal{R}$ is closed under the Souslin operation.}

\subsection{There are no hereditary productive $\gamma$-spaces}\label{FJ}
We show that if $X$ is an uncountable productive $\gamma$-set [F.~Jordan,
Productive local properties of function spaces,{\em Topology Appl.} {\bf 154}, 870--883, 2007], then there is a countable $Y\subseteq X$
such that $X\setminus Y$ is not Hurewicz.

Along the way we will prove a general result about Fr\'echet-$\alpha_2$ filters to gain information
about countable unions of $\gamma$-spaces and productive $\gamma$-spaces.
In particular, we answer a
question of A. Miller by showing that an increasing countable union of $\gamma$-spaces is again a $\gamma$-space.
We also use recent methods of B. Tsaban and L. Zdomskyy to show that $\lambda$-spaces with
Hurewicz property are precisely those spaces for which every co-countable set is Hurewicz.

\nby{Francis Jordan}

\section{Problem of the Issue}

For a sequence $\seq{X_n}$ of subsets of $X$, define
$\liminf X_n = \Union_m\allowbreak\bigcap_{n\ge m} X_n$.
For a family $\cF$ of subsets of $X$, $L(\cF)$ denotes
its closure under the operation $\liminf$.
$X$ has the \emph{$\delta$-property} if
for each open $\w$-cover $\cU$ of $X$, $X\in L(\cU)$.
This property was introduced by Gerlits and Nagy in their seminal paper \cite{GN}.

Clearly, the $\gamma$-property $\binom{\Omega}{\Gamma}$ implies the $\delta$-property.
$\sone(\Omega,\Gamma)=\binom{\Omega}{\Gamma}$ \cite{GN}.

\begin{prob}[Gerlits-Nagy \cite{GN}]\label{GNdelta}
Is the $\delta$-property equivalent to $\binom{\Omega}{\Gamma}$?
\end{prob}

Miller suggested that, as a union of an increasing sequence of
sets with the $\delta$-property has again the $\delta$-property,
one can obtain a negative answer by finding an increasing sequence
$\seq{X_n}$ of sets, each satisfying $\binom{\Omega}{\Gamma}$,
such that their union does not satisfy $\binom{\Omega}{\Gamma}$.
Recently, Francis Jordan solved Miller's question by proving that this is impossible
(Section \ref{FJ} above).

A closely-related problem, due to Sakai, concerns the Pytkeev property in function spaces.
This problem asks whether, for each set of reals $X$,
if $C_p(X)$ has the Pytkeev property then $C_p(X)$ is Fr\'echet (i.e., $X$ satisfies $\binom{\Omega}{\Gamma}$) -- see \cite{Sakai06}
for details.

\nby{Boaz Tsaban}

\ed